\documentclass[12pt,a4paper]{article}
\usepackage{amsmath,amsthm,amssymb, mathrsfs}
\usepackage{hyperref}
\usepackage{tikz}
\usepackage{graphicx}
\usepackage{cite}
\input xy
\xyoption{all}
\usepackage{epsfig}
\newtheorem{thm}{Theorem}[section]

 \theoremstyle{definition}

\theoremstyle{remark}

\numberwithin{equation}{section}
\newcommand{\ben}{\begin{enumerate}}

\newcommand{\een}{\end{enumerate}}
\newcommand{\bit}{\begin{itemize}}
\newcommand{\eit}{\end{itemize}}

\begin{document}

\title
{The Spherical Mean Transform with Data on a Parabola in the Plane}
\date{}

\vskip 1cm
\author{Yehonatan Salman \\ Email: salman.yehonatan@gmail.com\\ Weizmann Institute of Science}
\date{}

\maketitle

\begin{abstract}

In this paper we deal with the problem of recovering functions from their spherical mean transform $\mathcal{R}$, which integrates functions on circles in the plane, in case where the centers of the circles of integration are located on a parabola $\mathcal{P}$ while their radii can be chosen arbitrarily. Using our data, on the values of $\mathcal{R}$ on $\mathcal{P}$, we show how to extract its values in the exterior of $\mathcal{P}$ in case where the functions in question have compact support inside $\mathcal{P}$. Hence, one can use known inversion formulas for $\mathcal{R}$ in the exterior of $\mathcal{P}$ in order to obtain a reconstruction formula.

\end{abstract}

\section{Introduction and Mathematical Background}

\hskip0.6cm Denote by $\Bbb R$ and $\Bbb R^{+}$ respectively the real line and the ray $[0,\infty)$. Denote by $\Bbb R^{n}$ the $n$ dimensional Euclidean space. For a continuous function $f$, defined in $\Bbb R^{2}$, define the spherical mean transform $\mathcal{R}f$ of $f$ by
\vskip-0.2cm
$$\hskip-10cm\mathcal{R}f:\Bbb R^{2}\times\Bbb R^{+}\rightarrow\Bbb R,$$
$$\hskip-7.8cm(\mathcal{R}f)(x, r) = \int_{-\pi}^{\pi}f\left(x + re^{i\theta}\right)rd\theta.$$
That is, at each point $(x, r)$ the function $\mathcal{R}f$ evaluates the integral of $f$ on the circle with the center at $x$ and radius $r$.

The spherical mean transform became a major object of study in Integral Geometry in the last decades where various results concerning uniqueness, inversion and range theorems for this integral transform have been obtained. The inversion problem was studied in case where for a function $f$ in question, which is defined in $\Bbb R^{2}$ (or more generally in $\Bbb R^{n}$), the spherical mean transform $\mathcal{R}f$ is restricted to a cylindrical surface of the form $\Sigma\times\Bbb R^{+}$ where $\Sigma$ is an algebraic curve in $\Bbb R^{2}$ (or more generally an algebraic hypersurface in $\Bbb R^{n}$). That is, our data consists of all the integrals of the function $f$ on circles (or more generally hyperspheres) with centers on $\Sigma$, while no restriction is imposed on the set of radii, and our aim is to extract $f$ from this data.

In the last two decades inversion methods for the spherical mean transform have been obtained for various quadratic curves and hypersurfaces $\Sigma$ (see \cite{1, 5, 6, 7, 8, 9, 10, 11, 12, 13, 14, 15, 16, 17, 18}). In this paper we will concentrate on the case where $\Sigma$ is a parabola in $\Bbb R^{2}$ which is defined by the following parametrization
\vskip-0.2cm
$$\hskip-6cm\mathcal{P}_{\eta_{0}} = \left\{x = \frac{1}{2}\left(\xi^{2} - \eta_{0}^{2}\right), y = \xi\eta_{0}: \xi\in\Bbb R\right\}$$
where $\eta_{0}$ is a positive real number. Algebraically, $\mathcal{P}_{\eta_{0}}$ is given by $y^{2} - 2\eta_{0}^{2}x - \eta_{0}^{4} = 0$
and by using the translation $x = x' - \eta_{0}^{2} / 2$ it is an easy exercise to show that by varying $\eta_{0}$ we obtain the set of parabolas $x = cy^{2}, c > 0$. For a given $\eta_{0} > 0$ our data consists of the integrals of $f$ on circles with centers on $\mathcal{P}_{\eta_{0}}$ and arbitrary radii and our aim is to reconstruct $f$ from this data.

In order to obtain a reconstruction method we will use the same procedure that was used in \cite{17} where an inversion formula was found, for the case where the centers of the circles of integration are given on an ellipse, by using the expansion of the Bessel function $J_{0}$ in elliptical coordinates. In this paper we use a slight modification of the method introduced in \cite{17} where now we use the expansion of the modified Bessel function of the second kind of order zero $K_{0}$ in parabolic coordinates which was found by Cohl and Volkmer in \cite[Theorem 2.2]{2}. However, in this paper the reconstruction method is given more implicitly since we are not going to extract each function $f$ in question from its spherical mean transform $\mathcal{R}f$ explicitly. Instead, we will show that if the support of $f$ is inside $\mathcal{P}_{\eta_{0}}$ then one can extract the values of $\mathcal{R}f$ at the exterior of $\mathcal{P}_{\eta_{0}}$. That is, we can extract the following values
$$\hskip-0.5cm(\mathcal{R}f)(x, r), r\geq0, x\in\mathcal{P}_{\eta_{0}}^{ext} := \left\{x = \frac{1}{2}\left(\xi^{2} - \eta^{2}\right), y = \xi\eta:\xi\in\Bbb R, \eta\geq\eta_{0}\right\}.$$
Hence, by extracting the values of $\mathcal{R}f$ in the exterior of $\mathcal{P}_{\eta_{0}}$ we can use known inversion formulas, where we take circles with centers which are located on a curve $\Sigma'\subset\mathcal{P}_{\eta_{0}}^{ext}$, in order to reconstruct $f$. For example, we can take $\Sigma'$ to be any vertical line $x = x_{0}$ where $x_{0} < -\eta_{0}^{2} / 2$ and then, since each function in question is supported in the half plane $x\geq x_{0}$ (since $f$ is supported inside $\mathcal{P}_{\eta_{0}}$), it is known that in this case (i.e., when the centers of circles of integration are located on $x = x_{0}$) the function $f$ can be reconstructed (see for example \cite{5,11}).

For other inversion formulas of back-projection type where $\Sigma$ is a parabola in the plane (or more generally a paraboloid in $\Bbb R^{n}$) see \cite{9, 16}.

The main result of this paper is given in Theorem 2.1 in Sect. 2. Before formulating Theorem 2.1 we will introduce some notations and definitions which will be used in the proof of this theorem.\\

For every point $x\in\Bbb R^{2}$ denote by $x(\xi, \eta)$ its presentation in parabolic coordinates:
\vskip-0.2cm
$$\hskip-6.1cm x(\xi, \eta) = \left(\frac{1}{2}\left(\xi^{2} - \eta^{2}\right), \xi\eta\right), \xi\in\Bbb R, \eta\geq0$$
where the parameters $\xi$ and $\eta$ are uniquely determined by the point $x$. Denote by $\rho(\xi, \eta, \xi', \eta')$ the distance
\vskip-0.2cm
$$\hskip-7.3cm\rho(\xi, \eta, \xi', \eta') = \left|x(\xi, \eta) - x(\xi', \eta')\right|.$$
Let $K_{0}$ be the modified Bessel function of the second kind of order zero and let $H_{\nu}$, $\nu\in\Bbb R$ be the Hermite function of order $\nu$. For a non negative integer $n$ the function $H_{n}$ coincides with the Hermite polynomial
\vskip-0.2cm
$$\hskip-8.5cm H_{n}(x) = ( - 1)^{n}e^{x^{2}}\frac{d^{n}}{dx^{n}}e^{-x^{2}}.$$
For the expression of $H_{\nu}$ where $\nu$ is an arbitrary real number see \cite[Sect. 2]{2}.
In \cite[Theorem 2.2]{2} the following identity
\vskip-0.2cm
$$\hskip-5.5cm K_{0}(k\cdot\rho(\xi, \eta, \xi', \eta'))
 = \sqrt{\pi}e^{\frac{k}{2}\left((\eta')^{2} - (\xi')^{2} - \eta^{2} - \xi^{2}\right)}$$ \begin{equation}\hskip1cm\times\sum_{n = 0}^{\infty}\frac{(-i)^{n}}{2^{n - 1}n!}H_{n}\left(\sqrt{k}\xi\right)H_{- n - 1}\left(\sqrt{k}\eta\right)H_{n}\left(\sqrt{k}\xi'\right)H_{n}\left(i\sqrt{k}\eta'\right)\end{equation}
was proved in case where $0\leq \eta'\leq\eta$ and $k > 0$. For a function $F$, defined in $\Bbb R^{+}$, define the Mellin transform $\mathcal{M}F$ of $F$ by
$$\hskip-6.8cm(\mathcal{M}F)(s) = \int_{0}^{\infty}y^{s - 1}F(y)dy, \Re s > 0$$
where it should be noted that the above integral might not converge for every $\Re s > 0$. For the Mellin transform we have the following inversion and convolution formulas (see \cite{3}, Chapter 8.2 and 8.3):
\begin{equation}\hskip-3.85cm F(r) = \mathcal{M}^{-1}(\mathcal{M}F)(r) = \frac{1}{2\pi i}\int_{\varrho - i\infty}^{\varrho + i\infty}r^{-s}\mathcal{M}(F)(s)ds,\end{equation}
\begin{equation}\hskip-5.75cm\mathcal{M}(F_{1}\star F_{2})(s) = (\mathcal{M}F_{1})(s)(\mathcal{M}F_{2})(1 - s)\end{equation}
where the convolution $F_{1}\star F_{2}$ is defined by
\vskip-0.2cm
$$\hskip-6.7cm (F_{1}\star F_{2})(x) = \int_{0}^{\infty}F_{1}(xx')F_{2}(x')dx'.$$
In the Mellin inversion formula (1.2) the point $\sigma$ can be any point in $(0,\infty)$ for which the Mellin transform of $F$ decays to zero uniformly on a strip which contains the complex line $\sigma + it, t\in\Bbb R$.

Using the relationship between the Fourier and Mellin transforms (see \cite[Chap. 8, Sect. 2]{3}) it can be easily checked that formula (1.2) is valid for every continuous function $F$ with compact support while the Mellin convolution formula (1.3) is valid in any domain in $\Bbb C$ for which both $(\mathcal{M}F_{1})(\cdot)$ and $(\mathcal{M}F_{2})(1 - \cdot)$ exist.

\section{The Main Result}

\hskip0.6cm For a continuous function $f$, defined in $\Bbb R^{2}$, our aim is to extract $\mathcal{R}f$ in the solid cylinder $\mathcal{P}_{\eta_{0}}^{ext}\times\Bbb R^{+}$ from its restriction on $\mathcal{P}_{\eta_{0}}\times\Bbb R^{+}$. That is, our data consists of the following values
\vskip-0.2cm
$$\hskip-5cm F(\xi, r) = (\mathcal{R}f)(x(\xi, \eta_{0}), r), (\xi, r)\in\Bbb R\times\Bbb R^{+}$$
and our aim is to evaluate $(\mathcal{R}f)(x(\xi, \eta), r)$ for every $(\xi, \eta)\in\Bbb R\times[\eta_{0}, \infty)$ and $r\geq0$. For this we have the following result.

\begin{thm}

Let $f$ be a continuous function, defined in $\Bbb R^{2}$, which is compactly supported inside the parabola $\mathcal{P}_{\eta_{0}}$, i.e., $f$ vanishes outside a bounded set in $\Bbb R^{2}$ and
$$\hskip-8.8cm f(x(\xi, \eta)) = 0, \eta\geq\eta_{0}.$$
Define
$$\hskip-1.6cm \Psi(\xi, \eta, k) = e^{-\frac{k}{2}\left(\eta^{2} + \xi^{2}\right)}\sum_{m = 0}^{\infty}H_{m}\left(\sqrt{k}\xi\right)H_{- m - 1}\left(\sqrt{k}\eta\right)(\Lambda_{m}f)(k)$$
for $(\xi, \eta)\in\Bbb R\times[\eta_{0}, \infty), k > 0$ where
\vskip-0.2cm
$$\hskip-6.25cm(\Lambda_{m}f)(k) = \frac{\sqrt{k}e^{\frac{k}{2}\eta_{0}^{2}}}{2^{m}m!\sqrt{\pi}H_{- m - 1}\left(\sqrt{k}\eta_{0}\right)}$$
$$\hskip1.5cm\times\int_{-\infty}^{\infty}\int_{0}^{\infty}(\mathcal{R}f)(x(\xi, \eta_{0}), r)K_{0}(kr)H_{m}\left(\sqrt{k}\xi\right)e^{-\frac{1}{2}\xi^{2}}drd\xi.$$
Then, for every $(\xi,\eta)\in\Bbb R\times[\eta_{0},\infty)$ and $r\geq0$ we have
\vskip-0.2cm
$$\hskip-3.1cm(\mathcal{R}f)(x(\xi, \eta), r) = \frac{1}{i\pi}\int_{\sigma - i\infty}^{\sigma + i\infty}\frac{2^{s}r^{-s}(\mathcal{M}\Psi)(\xi, \eta, 1 - s)ds}{\Gamma^{2}((1 - s) / 2)}$$
where $\sigma$ can be any point in the interval $(0,1)$ and where $\Gamma$ denotes the Gamma function.

\end{thm}

\vskip0.2cm

\begin{proof}

Using the definition of the spherical mean transform for every $\xi\in\Bbb R$ and $k > 0$ we have
$$\hskip-0.5cm\int_{0}^{\infty}(\mathcal{R}f)(x(\xi, \eta_{0}), r)K_{0}(kr)dr = \int_{0}^{\infty}\int_{-\pi}^{\pi}f\left(x(\xi, \eta_{0}) + re^{i\theta}\right)rd\theta K_{0}(kr)dr $$
$$\hskip-1.5cm = [y = x(\xi, \eta_{0}) + re^{i\theta}, dy = rd\theta dr] = \int_{\Bbb R^{2}}f(y)K_{0}\left(k\left|y - x(\xi, \eta_{0})\right|\right)dy$$
$$\hskip-6.05cm = \left[y = y(\xi', \eta'), dy = ((\xi')^{2} + (\eta')^{2})d\xi'd\eta'\right]$$
\begin{equation}\hskip-5.2cm = \int_{0}^{\eta_{0}}\int_{-\infty}^{\infty}f^{\ast}(\xi', \eta')K_{0}\left(k\cdot\rho(\xi, \eta_{0}, \xi', \eta')\right)d\xi'd\eta'\end{equation}
where we denote
\vskip-0.2cm
$$\hskip-7cm f^{\ast}(\xi', \eta') = ((\xi')^{2} + (\eta')^{2})f(y(\xi', \eta'))$$
and in the last passage of equation (2.1) we used the fact that $f$ is supported inside $\mathcal{P}_{\eta_{0}}$. Using the expansion (1.1) of $K_{0}$ and the fact that $\eta'\leq\eta_{0}$ we have
$$\hskip-7.8cm\int_{0}^{\infty}(\mathcal{R}f)(x(\xi, \eta_{0}), r)K_{0}(kr)dr$$
$$\hskip-3.5cm = \sqrt{\pi}e^{-\frac{k}{2}(\eta_{0}^{2} + \xi^{2})}\sum_{n = 0}^{\infty}\frac{(- i)^{n}}{2^{n - 1}n!}H_{n}\left(\sqrt{k}\xi\right)H_{- n - 1}\left(\sqrt{k}\eta_{0}\right)$$
$$\hskip-1.5cm\times\int_{0}^{\eta_{0}}\int_{-\infty}^{\infty}f^{\ast}(\xi', \eta')e^{\frac{k}{2}\left((\eta')^{2} - (\xi')^{2}\right)}H_{n}\left(\sqrt{k}\xi'\right)H_{n}\left(i\sqrt{k}\eta'\right)d\xi'd\eta'.$$
Using the following orthogonality relations
$$\hskip-3.25cm\int_{-\infty}^{\infty}H_{m}\left(\sqrt{k}x\right)H_{n}\left(\sqrt{k}x\right)e^{-kx^{2}}dx = \sqrt{\pi / k}2^{n}n!\delta_{n,m}$$
for the Hermite polynomials we obtain that
$$\hskip-3.5cm\int_{-\infty}^{\infty}\int_{0}^{\infty}(\mathcal{R}f)(x(\xi, \eta_{0}), r)K_{0}(kr)H_{m}\left(\sqrt{k}\xi\right)e^{-\frac{k}{2}\xi^{2}}drd\xi$$
$$ \hskip-7.2cm = \frac{2(- i)^{m}\pi}{\sqrt{k}}\cdot e^{-\frac{k}{2}\eta_{0}^{2}}H_{- m - 1}\left(\sqrt{k}\eta_{0}\right)$$
\begin{equation}\hskip-2cm\times\int_{0}^{\eta_{0}}\int_{-\infty}^{\infty}f^{\ast}(\xi', \eta')e^{\frac{k}{2}\left((\eta')^{2} - (\xi')^{2}\right)}H_{m}\left(\sqrt{k}\xi'\right)H_{m}\left(i\sqrt{k}\eta'\right)d\xi'd\eta'.\end{equation}
Let us denote
\vskip-0.2cm
$$\hskip-6.7cm(\Lambda_{m}f)(k) = \frac{\sqrt{k}e^{\frac{k}{2}\eta_{0}^{2}}}{2^{m}m!\sqrt{\pi}H_{- m - 1}\left(\sqrt{k}\eta_{0}\right)}$$ $$\hskip-0.2cm\times\int_{-\infty}^{\infty}\int_{0}^{\infty}(\mathcal{R}f)(x(\xi, \eta_{0}), r)K_{0}(kr)H_{m}\left(\sqrt{k}\xi\right)e^{-\frac{1}{2}\xi^{2}}drd\xi$$
then from equation (2.2) for every $\xi''\in\Bbb R$ and $\eta''\geq\eta_{0}$ we have
\vskip-0.2cm
$$\hskip-3cm e^{-\frac{k}{2}\left((\eta'')^{2} + (\xi'')^{2}\right)}\sum_{m = 0}^{\infty}H_{m}\left(\sqrt{k}\xi''\right)H_{- m - 1}\left(\sqrt{k}\eta''\right)(\Lambda_{m}f)(k)$$
$$\hskip-2.5cm = \sqrt{\pi}e^{-\frac{k}{2}((\eta'')^{2} + (\xi'')^{2})}\sum_{m = 0}^{\infty}\frac{(- i)^{m}}{2^{m - 1}m!}H_{m}\left(\sqrt{k}\xi''\right)H_{- m - 1}\left(\sqrt{k}\eta''\right)$$
$$\hskip-2cm\times\int_{0}^{\eta_{0}}\int_{-\infty}^{\infty}f^{\ast}(\xi', \eta')e^{\frac{k}{2}\left((\eta')^{2} - (\xi')^{2}\right)}H_{m}\left(\sqrt{k}\xi'\right)H_{m}\left(i\sqrt{k}\eta'\right)d\xi'd\eta'$$
$$\hskip-4.95cm = \int_{0}^{\eta_{0}}\int_{-\infty}^{\infty}\sqrt{\pi}f^{\ast}(\xi', \eta')e^{\frac{k}{2}\left((\eta')^{2} - (\xi')^{2} - (\eta'')^{2} - (\xi'')^{2}\right)}$$ $$\hskip0.5cm\times\left[\sum_{m = 0}^{\infty}\frac{(- i)^{m}}{2^{m - 1}m!}H_{m}\left(\sqrt{k}\xi''\right)H_{- m - 1}\left(\sqrt{k}\eta''\right)H_{m}\left(\sqrt{k}\xi'\right)H_{m}\left(i\sqrt{k}\eta'\right)\right]d\xi'd\eta'$$
$$\hskip-4.7cm \underset{(*)}{=} \int_{0}^{\eta_{0}}\int_{-\infty}^{\infty}K_{0}(k\cdot\rho(\xi'', \eta'', \xi', \eta'))f^{\ast}(\xi', \eta')d\xi'd\eta' $$
$$ \hskip-5.7cm = \left[y = y(\xi', \eta'), dy = ((\xi')^{2} + (\eta')^{2})d\xi'd\eta'\right]$$
\begin{equation}\hskip-7.25cm = \int_{\Bbb R^{2}}K_{0}\left(k\left|x(\xi'', \eta'') - y\right|\right)f(y)dy\end{equation}
where in the passage $(\ast)$ we used the expansion (1.1) of the modified Bessel function of the second kind $K_{0}$ and that $\eta'\leq\eta''$ (since $\eta'\leq\eta_{0}\leq\eta''$). Making the change of variables
\vskip-0.2cm
$$\hskip-7.2cm y = x(\xi'', \eta'') + re^{i\theta}, dy = rd\theta dr$$
in equation (2.3) we obtain
$$\hskip-3cm e^{-\frac{k}{2}\left((\eta'')^{2} + (\xi'')^{2}\right)}\sum_{m = 0}^{\infty}H_{m}\left(\sqrt{k}\xi''\right)H_{- m - 1}\left(\sqrt{k}\eta''\right)(\Lambda_{m}f)(k)$$
$$\hskip-5.5cm = \int_{0}^{\infty}\int_{-\pi}^{\pi}K_{0}(kr)f\left(x(\xi'', \eta'') + re^{i\theta}\right)rd\theta dr$$
\begin{equation}\hskip-7.15cm = \int_{0}^{\infty}K_{0}(kr)(\mathcal{R}f)(x(\xi'', \eta''), r)dr.\end{equation}
If we denote the function in the left hand side of equation (2.4) by $\Psi = \Psi(\xi'', \eta'', k)$ then by taking the Mellin transform on both sides of (2.4) with respect to the variable $k$, using the Melling convolution formula (1.3) and using the fact that the Mellin transform of $K_{0}$ is given by
\vskip-0.2cm
$$\hskip-7cm(\mathcal{M}K_{0})(s) = 2^{s - 2}\Gamma^{2}\left(\frac{s}{2}\right), \Re s > 0$$
where $\Gamma$ denotes the Gamma function (see \cite[Chap. VI, Sect 6.8, formula 26]{4}) we obtain
\vskip-0.2cm
$$\hskip-3.1cm(\mathcal{M}\Psi)(\xi'', \eta'', s) = 2^{s - 2}\Gamma^{2}(s / 2)\mathcal{M}\left(\mathcal{R}f\right)(x(\xi'', \eta''), 1 - s)$$
for $0 < \Re s < 1$. Equivalently, we have
$$\hskip-4.2cm\mathcal{M}\left(\mathcal{R}f\right)(x(\xi'', \eta''), s) = \frac{2^{s + 1}(\mathcal{M}\Psi)(\xi'', \eta'', 1 - s)}{\Gamma^{2}((1 - s) / 2)}$$
where the last equality is true for $|\xi''| < \infty, \eta''\geq\eta_{0}$ and $0 < \Re s < 1$. For every fixed $x(\xi'', \eta'')$, where $|\xi''| < \infty, \eta''\geq\eta_{0}$, the function $(\mathcal{R}f)(x(\xi'', \eta''), \cdot)$ is continuous and has compact support (since $f$ is continuous and has compact support). Hence, we can apply the inverse Mellin transform and thus Theorem 2.1 is proved.

\end{proof}


\begin{thebibliography}{}

\bibitem{1} Y. A. Antipov, R. Estrada and B. Rubin. \textit{Inversion formulas for spherical means in constant curvature spaces}. Journal D'Anal. Math.
    118, 623-656, 2012.

\bibitem{2} H.S. Cohl and H.Volkmer. \textit{Eigenfunction expansions for a fundamental solution of Laplace's equation on $\mathbf{R^{3}}$ in parabolic and elliptic cylinder coordinates}. J. Phys. A: Math. Theor. 45 355204, 2012.

\bibitem{3} B. Dambaru and L. Debnath. \textit{Integral Transforms and Their Applications}, CRC Press, New York, 2007.

\bibitem{4} A. Erd\`{e}lyi A, W. Magnus, F. Oberhettinger F and F.G. Tricomi. \textit{Tables of integral transforms}, New
York: McGraw-Hill; 1954.

\bibitem{5} A. J. Fawcett. \textit{Inversion of N-Dimensional Spherical Averages}, SIAM J. Appl. Math. 45(2), 1985.
336-341

\bibitem{6} D. Finch, M. Haltmeier, and Rakesh. \textit{Inversion of spherical means and the wave equation in
              even dimensions}. SIAM J. Appl. Math. 68, no.2, 392-412, 2007.

\bibitem{7} D. Finch, S. Patch, and Rakesh. \textit{Determining a function from its mean values over a family of spheres}. SIAM J. Math. Anal, 35, no.
    5, 1213-1240, 2004.

\bibitem{8} M. Haltmeier. \textit{Exact reconstruction formula for the spherical mean Radon transform on ellipsoids}. Inverse Probl. 30,
           no. 10, 105006, 2014.

\bibitem{9} M. Haltmeier and S. Pereverzyev Jr. \textit{Recovering a function from circular means or wave data on the boundary of parabolic domains}. SIAM J. Imaging Sci., 8(1):592-610, 2015.

\bibitem{10} L. A. Kunyansky. \textit{Explicit inversion formulae for the spherical mean Radon transform}. Inverse Probl. 23, no. 1, 373-383, 2007.

\bibitem{11} E. k. Narayanan and Rakesh. \textit{Spherical means with centers on a hyperplane in even dimensions},
Inverse Problems 26(3) 035014, 2010.

\bibitem{12} F. Natterer. \textit{Photo-acoustic inversion in convex domains}. Inverse Probl. Imaging, no. 2, 315-320, 2012.

\bibitem{13} L. V. Nguyen. \textit{A family of inversion formulas for thermoacoustic tomography}. Inverse Probl. Imaging, no. 4, 649-675, 2009.

\bibitem{14} S. J. Norton. \textit{Reconstruction of a two-dimensional reflecting medium over a circular domain: exact solution}. J. Acoust. Soc. Amer.
    67, 1266-1273, 1980.

\bibitem{15} V. P. Palamodov. \textit{A uniform reconstruction formula in integral geometry}. Inverse Probl.,
           no. 6, 065014, 2012.

\bibitem{16} V. P. Palamodov. \textit{Reconstruction from Integral Data}. Monographs and Research Notes in
Mathematics. CRC Press, Boca Raton, 2016.

\bibitem{17} Y. Salman. \textit{Recovering Functions from the Spherical Mean Transform with Data on an Ellipse Using Eigenfunction Expansion in Elliptical Coordinates}. Anal. Math. Phys., Advance online publication, 2017.

\bibitem{18} M. Xu and L. V. Wang. \textit{Universal back-projection algorithm for photoacoustic computed
tomography}. Physical Review E, 71, 2005.

\end{thebibliography}
\end{document}